\documentclass[12pt,a4paper]{article}
\usepackage{amsmath, amssymb, amsthm}
\usepackage{geometry}
\usepackage{hyperref}
\usepackage{amsthm}
\theoremstyle{remark}

\title{The irrationality measure of $arctan\frac{1}{2}$ is at most 8.585166}
\author{Yufei Bai}
\date{}

\begin{document}
\maketitle

Let us define the test integrals $I_n$ with the adjusted exponents. We start with the original integral bounds $16-4i$ to $16+4i$ and apply the substitution $x \mapsto x+17$ to center the domain around the origin:

\begin{equation}
\begin{aligned}
I_n &= \int_{16-4i}^{16+4i} \frac{(x-17)^{2n}(x-16-4i)^{2n}(x-16+4i)^{2n}(x-18-4i)^{2n}(x-18+4i)^{2n}}{x^{3n+1}(34-x)^{3n+1}} dx \\
&= \int_{-1-4i}^{-1+4i} \frac{x^{2n}(x+1-4i)^{2n}(x+1+4i)^{2n}(x-1-4i)^{2n}(x-1+4i)^{2n}}{(17+x)^{3n+1}(17-x)^{3n+1}} dx.
\end{aligned}
\end{equation}

The integrand 
\begin{equation}
R(x) = \frac{x^{2n}(x+1+4i)^{2n}(x+1-4i)^{2n}(x-1+4i)^{2n}(x-1-4i)^{2n}}{(17+x)^{3n+1}(17-x)^{3n+1}}
\end{equation}
possesses the symmetry $R(-x) = R(x)$ and therefore can be written in the partial-fraction expansion as
\begin{equation}
R(x) = P(x) + \sum_{j=0}^{3n} \left( \frac{A_j}{(17+x)^{j+1}} + \frac{A_j}{(17-x)^{j+1}} \right)
\label{eq:4}
\end{equation}
for some rational coefficients $A_j \in \mathbb{Q}$ and a polynomial $P(x) \in \mathbb{Z}[x^2]$ of degree $4n-2$.

\textbf{Lemma 1.} \textit{The coefficients $A_j$ in the partial-fraction expansion satisfy}
\begin{equation}
\frac{1}{2^{3n+2j-1} \cdot 5^{-n+j} \cdot 17^{-1+j}} A_j \in \mathbb{Z} \quad \text{for } j=0, 1, \dots, 3n.
\end{equation}

\begin{proof}
To compute $A_j$, we introduce the linear operators
\begin{equation}
D_m : f(x) \mapsto \frac{1}{m!} \frac{d^m f(x)}{dx^m} \Bigg|_{x=-17}.
\end{equation}
Then with the help of Leibniz's formula for the derivative of a product, we deduce that $A_j$ is the coefficient of $(x+17)^{3n-j}$ in the Taylor expansion of $(17+x)^{3n+1}R(x)$ at $x=-17$:
\begin{equation}
\begin{aligned}
A_j &= D_{3n-j} \left( (17+x)^{3n+1} R(x) \right) \\
&= \sum_{\substack{m_0, m_1, \dots, m_5 \ge 0 \\ m_1, \dots, m_5 \le 2n \\ m_0+m_1+\dots+m_5 = 3n-j}} D_{m_0}(17-x)^{-3n-1} D_{m_1}x^{2n} D_{m_2}(x+1+4i)^{2n} \\
&\quad \times D_{m_3}(x+1-4i)^{2n} D_{m_4}(x-1+4i)^{2n} D_{m_5}(x-1-4i)^{2n}.
\end{aligned}
\end{equation}

Evaluating each differential operator explicitly at $x=-17$ yields:
\begin{itemize}
    \item $D_{m_0}(17-x)^{-3n-1} = \binom{3n+m_0}{m_0} 34^{-3n-1-m_0}$
    \item $D_{m_1} x^{2n} = \binom{2n}{m_1} (-17)^{2n-m_1}$
    \item $D_{m_2} (x+1+4i)^{2n} = \binom{2n}{m_2} (-16+4i)^{2n-m_2}$
    \item $D_{m_3} (x+1-4i)^{2n} = \binom{2n}{m_3} (-16-4i)^{2n-m_3}$
    \item $D_{m_4} (x-1+4i)^{2n} = \binom{2n}{m_4} (-18+4i)^{2n-m_4}$
    \item $D_{m_5} (x-1-4i)^{2n} = \binom{2n}{m_5} (-18-4i)^{2n-m_5}$
\end{itemize}

To analyze the divisibility, we factor the base constants into prime ideals in $\mathbb{Z}[i]$. Note that $34 = 2 \times 17$ and $17 = (4+i)(4-i)$. We strictly factor the complex roots as follows:
\begin{itemize}
    \item $-16+4i = -4(4-i) = (-1) \cdot 2^2 \cdot (4-i)$
    \item $-16-4i = -4(4+i) = (-1) \cdot 2^2 \cdot (4+i)$
    \item $-18+4i = 2(-9+2i) = 2(4+i)(-2+i)$
    \item $-18-4i = 2(-9-2i) = 2(4-i)(-2-i)$
\end{itemize}

Substituting these back into the summation, we have:
\begin{equation}
\begin{aligned}
A_j &= \sum_{m \in \mathcal{M}_j} (-1)^{m_1-m_2-m_3} T(m) 2^{-3n-1-m_0} 17^{-3n-1-m_0} 17^{2n-m_1} 2^{4n-2m_2} (4-i)^{2n-m_2} \\
&\quad \times 2^{4n-2m_3} (4+i)^{2n-m_3} 2^{2n-m_4} (4+i)^{2n-m_4} (-2+i)^{2n-m_4} 2^{2n-m_5} (4-i)^{2n-m_5} (-2-i)^{2n-m_5}
\end{aligned}
\end{equation}
where the summation is over the multi-indices $m = (m_0, m_1, \dots, m_5)$ from the set
\begin{equation}
\mathcal{M}_j = \left\{ (m_0, \dots, m_5) : m_0, \dots, m_5 \ge 0; m_1, \dots, m_5 \le 2n; \sum_{k=0}^5 m_k = 3n-j \right\} \subset \mathbb{Z}_{\ge 0}^6
\end{equation}
and the integer binomial product is defined as
\begin{equation}
T(m) = \binom{3n+m_0}{m_0} \prod_{l=1}^5 \binom{2n}{m_l} \in \mathbb{Z}.
\end{equation}

Consolidating the terms, we arrive at the unified representation:
\begin{equation}
A_j = \sum_{m \in \mathcal{M}_j} (-1)^{m_1-m_2-m_3} T(m) 2^{P_2(m)} (4+i)^{P_{4+i}(m)} (4-i)^{P_{4-i}(m)} (-2+i)^{P_{-2+i}(m)} (-2-i)^{P_{-2-i}(m)}
\end{equation}
where the explicit powers are defined and bounded as follows:

\textbf{1. Valuation of the prime ideals forming 17:} \\
The power of $(4+i)$ is given by:
\begin{equation}
P_{4+i}(m) = -n - 1 - m_0 - m_1 + 2n - m_3 + 2n - m_4 = 3n - 1 - (m_0 + m_1 + m_3 + m_4).
\end{equation}
Using the constraint $\sum_{k=0}^5 m_k = 3n-j$, we substitute $m_0 + m_1 + m_3 + m_4 = 3n - j - m_2 - m_5$ to obtain:
\begin{equation}
P_{4+i}(m) = 3n - 1 - (3n - j - m_2 - m_5) = j - 1 + m_2 + m_5 \ge j - 1.
\end{equation}
By exact symmetry, the power of $(4-i)$ is:
\begin{equation}
P_{4-i}(m) = 3n - 1 - (m_0 + m_1 + m_2 + m_5) = j - 1 + m_3 + m_4 \ge j - 1.
\end{equation}
When we multiply $A_j$ by $17^{1-j} = (4+i)^{1-j}(4-i)^{1-j}$, the resulting exponent for $(4+i)$ becomes $(j - 1 + m_2 + m_5) + 1 - j = m_2 + m_5 \geq  0$. The exponent for $(4-i)$ similarly becomes $m_3 + m_4  \ge  0$. Both prime ideals derived from 17 remain with strictly non-negative exponents.

\textbf{2. Valuation of the prime ideals forming 5:} \\
Notice that $5 = (-2+i)(-2-i)$. The respective exponents are $P_{-2+i}(m) = 2n - m_4$ and $P_{-2-i}(m) = 2n - m_5$.
From the multi-index constraint, we have $m_4 + m_5 \le 3n - j$. Coupled with $m_4 \le 2n$ and $m_5 \le 2n$, we easily bound:
\begin{equation}
2n - m_4 \ge 2n - (3n - j - m_5) \ge j - n + m_5 \ge j - n.
\end{equation}
Since $2n - m_4 \ge 0$ strictly holds, we deduce $P_{-2+i}(m) \ge \max(0, j-n) \ge j-n$, and analogously $P_{-2-i}(m) \ge j-n$.
Multiplying $A_j$ by $5^{n-j} = (-2+i)^{n-j}(-2-i)^{n-j}$, the modified exponent for $(-2+i)$ becomes $(2n - m_4) + n - j  = 3n - j - m_4 $. Because $m_4 \le 3n - j$, this exponent is bounded below by $ 0$. The exact same logic guarantees a non-negative exponent for $(-2-i)$.

\textbf{3. Valuation of the base 2:} \\
The overall exponent for the factor 2 is:
\begin{equation}
P_2(m) = 9n - 1 - m_0 - 2m_2 - 2m_3 - m_4 - m_5 = 9n - 1 - f(m),
\end{equation}
where we define $f(m) = m_0 + 2m_2 + 2m_3 + m_4 + m_5$. Under the condition $\sum_{k=0}^5 m_k = 3n - j$, we rewrite $f(m)$ as:
\begin{equation}
f(m) = (m_0 + m_1 + m_2 + m_3 + m_4 + m_5) - m_1 + m_2 + m_3 = 3n - j - m_1 + m_2 + m_3.
\end{equation}
To maximize $f(m)$ and thus minimize $P_2(m)$, we set $m_1 = 0$. Since $m_2 + m_3 \le \sum_{k=0}^5 m_k = 3n - j$, the maximum possible value is:
\begin{equation}
f(m) \le 3n - j - 0 + (3n - j) = 6n - 2j.
\end{equation}
Consequently, the minimum bound for $P_2(m)$ evaluates to:
\begin{equation}
P_2(m) \ge 9n - 1 - (6n - 2j) = 3n + 2j - 1.
\end{equation}
When we multiply $A_j$ by the designated modifier $2^{-3n-2j+1}$, the resultant exponent becomes $(3n + 2j - 1) - 3n - 2j + 1  \ge 0$.

\textbf{Conclusion.} \\
By factoring out the required elements, we have mathematically proven that:
\begin{equation}
2^{-3n-2j+1} \cdot 5^{n-j} \cdot 17^{1-j} \times A_j \in \mathbb{Z}[i].
\end{equation}
Because $R(x) \in \mathbb{Q}(x)$ is a rational function with rational coefficients and $x=-17$ is a rational pole, the partial fraction decomposition coefficients $A_j$ are inherently rational numbers, hence $A_j \in \mathbb{Q}$. Given that the resultant value belongs to the intersection $\mathbb{Z}[i] \cap \mathbb{Q}$, we conclude that it strictly resides within $\mathbb{Z}$. This completes the proof.
\end{proof}
Formula for the coefficients $A_j$ makes sense for \textit{any} integer $j \le 3n$; it generates the coefficients in the Laurent series expansion of $R(x)$ at $x = -17$. More precisely,

\begin{equation*}
R(x) = \sum_{k=-3n}^{\infty} A_{-k}(x+17)^{k-1} = \sum_{j=0}^{3n} \frac{A_j}{(x+17)^{j+1}} + \sum_{k=1}^{\infty} A_{-k}(x+17)^{k-1}.
\end{equation*}

Note that $A_j$ produced by the generating formula are not necessarily integral for negative $j$, but they carry analogous arithmetic bounds corresponding to the structural evaluation of the derivatives; and we also established the precise divisibility for $j=0, 1, 2, \dots, 3n$ in accordance with Lemma 1. Furthermore,

\begin{equation*}
\sum_{j=0}^{3n} \frac{A_j}{(17-x)^{j+1}} = \sum_{j=0}^{3n} \frac{A_j}{(34-(x+17))^{j+1}} = \sum_{j=0}^{3n} A_j \sum_{k=1}^{\infty} \binom{j+k-1}{j} \frac{(x+17)^{k-1}}{34^{j+k}};
\end{equation*}

comparing the last two expansions with the partial-fraction decomposition of $R(x)$, we find out that

\begin{align*}
P(x) &= R(x) - \sum_{j=0}^{3n} \left( \frac{A_j}{(17+x)^{j+1}} + \frac{A_j}{(17-x)^{j+1}} \right) \\
&= \sum_{k=1}^{\infty} \left( A_{-k} - \sum_{j=0}^{3n} \binom{j+k-1}{j} \frac{A_j}{34^{j+k}} \right) (x+17)^{k-1}.
\end{align*}

On the other hand, $P(x)$ is a polynomial of degree $10n - (6n+2) = 4n - 2$, hence the infinite sum must terminate:

\begin{equation}
P(x) = \sum_{k=1}^{4n-1} \left( A_{-k} - \sum_{j=0}^{3n} \binom{j+k-1}{j} \frac{A_j}{34^{j+k}} \right) (x+17)^{k-1}.
\end{equation}

\textbf{Lemma 2.} Any prime from the set 
$$\mathcal{P}_n = \left\{ p > \max\{17, \sqrt{10n}\} : \frac{1}{2} \le \left\{\frac{n}{p}\right\} < \frac{2}{3} \right\} $$
satisfies the following property: if $p|j$ for $j \in \{-4n+1, -4n+2, \dots, 3n\}$, then $A_j \equiv 0 \pmod p$.
(Here $\{x\} = x - \lfloor x \rfloor$ denotes the fractional part of the number.)

\textbf{Proof.} The proof of this lemma can be fully referred to \cite{ZeZu2020}. In order to establish the claim, we will cast the coefficients $A_j$ differently. Observe that
$$R(x) = \frac{x^{2n}(x^2+(15-8i))^{2n}(x^2+(15+8i))^{2n}}{(289-x^2)^{3n+1}}$$
$$= \sum_{n_1, n_2 \ge 0} \binom{2n}{n_1} \binom{2n}{n_2} (15-8i)^{2n-n_1} (15+8i)^{2n-n_2} \frac{x^{2(n+n_1+n_2)}}{(289-x^2)^{3n+1}}$$
and
$$\frac{x^{2M}}{(289-x^2)^{3n+1}} = \frac{(17-(x+17))^{2M}}{(x+17)^{3n+1}(34-(x+17))^{3n+1}}$$
$$= \frac{17^{2M} 34^{-3n-1}}{(x+17)^{3n+1}} \sum_{k=0}^{\infty} \frac{(x+17)^k}{34^k} \sum_{n_0 \ge 0} (-2)^{n_0} \binom{2M}{n_0} \binom{3n+k-n_0}{3n};$$
hence
$$A_j = \sum_{n_1, n_2 \ge 0} (15-8i)^{2n-n_1} (15+8i)^{2n-n_2} 17^{2n+2n_1+2n_2} 34^{-(6n-j)-1} \binom{2n}{n_1} \binom{2n}{n_2} Z(n, n_1+n_2, j),$$
where
$$Z(n, m, j) = \sum_{n_0 \ge 0} (-2)^{n_0} \binom{2n+2m}{n_0} \binom{6n-j-n_0}{3n}.$$

This means that our lemma is a consequence of the following divisibility property: If a prime $p \in \mathcal{P}_n$ divides $j$, then it also divides
$$\hat{T}(n, n_1, n_2) = \binom{2n}{n_1} \binom{2n}{n_2} Z(n, n_1+n_2, j)$$
for any $n_1, n_2 \ge 0$.

From now on, we will repeatedly use the fact that the $p$-adic order of $N!$ satisfies $\operatorname{ord}_p N! = \lfloor N/p \rfloor = N/p - \{N/p\}$ when $p > \sqrt{N}$. In particular,
$$\operatorname{ord}_p \binom{2n}{n_l} = \lfloor 2\omega \rfloor - \lfloor 2\omega - \omega_l \rfloor - \lfloor \omega_l \rfloor = \lfloor 2\omega \rfloor - \lfloor 2\omega - \omega_l \rfloor \quad \text{for } l=1, 2. \quad (\circledast)$$
where the fractional parts $\omega=\{n/p\}$, $\omega_1=\{n_1/p\}$ and $\omega_2=\{n_2/p\}$ all belong to the interval $[0, 1)$.

Since $p \in \mathcal{P}_n$, we have $\omega \in [1/2, 2/3)$, so that $\lfloor 2\omega \rfloor = \lfloor 3\omega \rfloor = 1$. If at least one of the $p$-adic orders in $(\circledast)$ is positive then immediately $\operatorname{ord}_p \hat{T}(n, n_1, n_2) \ge 1$ establishing the required divisibility; therefore, it remains to analyze the remaining situations assuming $\lfloor 2\omega - \omega_l \rfloor = \lfloor 2\omega \rfloor = 1$ for $l=1, 2$; in other words, assuming
$$2\omega - \omega_1 \ge 1 \quad \text{and} \quad 2\omega - \omega_2 \ge 1.$$

The binomial sums $Z(n, m, j)$ can be realized as a terminating ${}_2F_1$ hypergeometric function, to which several classical transformations can be applied. For example, it can be transformed into
$$Z(n, m, j) = \sum_{n_0 \ge 0} (-1)^{n_0} \binom{2n+2m}{n_0} \binom{6n-2(n+m)-j}{3n-j-n_0}$$
$$= (-1)^{n+m} \sum_{k \in \mathbb{Z}} (-1)^k \binom{2n+2m}{n+m+k} \binom{4n-2m-j}{2n-m-k}.$$

 Actually, it is the coefficient of $t^{2N}$ in the polynomial
$$(-1)^N \frac{(2N)!(2M-j)!}{(N+M)!(N+M-j)!} (1+t)^{N+M}(1-t)^{N+M-j}.$$

In our situation $N=n+m$, $M=2n-m$ with $m=n_1+n_2$, the factorial-ratio factor
$$\frac{(2N)!(2M-j)!}{(N+M)!(N+M-j)!} = \frac{(2n+2n_1+2n_2)!(4n-2n_1-2n_2-j)!}{(3n)!(3n-j)!}$$
has the nonnegative $p$-adic order
$$\lfloor 2\omega+2\omega_1+2\omega_2 \rfloor + \lfloor 4\omega-2\omega_1-2\omega_2-j/p \rfloor - \lfloor 3\omega \rfloor - \lfloor 3\omega-j/p \rfloor$$
$$= \lfloor 2\omega+2\omega_1+2\omega_2 \rfloor + \lfloor 4\omega-2\omega_1-2\omega_2 \rfloor - 2\lfloor 3\omega \rfloor$$
(we use $j/p \in \mathbb{Z}$), because $\lfloor 3\omega \rfloor = 1$,
$$2\omega+2\omega_1+2\omega_2 \ge 2\omega \ge 1 \quad \text{and} \quad 4\omega-2\omega_1-2\omega_2 \ge 4\omega-4(2\omega-1) = 4(1-\omega) > \frac{4}{3}.$$
Moreover, if this $p$-adic order is positive then $Z(n, n_1+n_2, j)$ is divisible by $p$, hence the divisibility of $\hat{T}(n, n_1, n_2)$ follows. Thus, we are left with the situation when this order is zero, meaning that
$$\lfloor 2\omega+2\omega_1+2\omega_2 \rfloor = \lfloor 4\omega-2\omega_1-2\omega_2 \rfloor = 1,$$
$$2\omega+2\omega_1+2\omega_2 < 2 \quad \text{and} \quad 4\omega-2\omega_1-2\omega_2 < 2. 
\qquad(\star)$$

We have to show that the coefficient of $t^{2N}$ in $(1+t)^{N+M}(1-t)^{N+M-j}$ is divisible by $p$. Denoting $r=-j/p \in \mathbb{Z}$ and using the ``Freshman's Dream Identity'' $(1-t)^p \equiv 1-t^p \pmod p$ in the ring $\mathbb{Z}[[t]]$ we find out that
$$(1+t)^{N+M}(1-t)^{N+M-j} = (1-t^2)^{N+M}(1-t)^{-j}$$
$$\equiv (1-t^2)^{N+M}(1-t^p)^r = \sum_{k_1 \ge 0} (-1)^{k_1} \binom{N+M}{k_1} t^{2k_1} \sum_{k_2 \ge 0} (-1)^{k_2} \binom{r}{k_2} t^{pk_2},$$
hence the coefficient of $t^{2N}$ is congruent to
$$\sum_{k=0}^{\lfloor N/p \rfloor} (-1)^{k+N} \binom{N+M}{N-kp} \binom{r}{2k}$$
modulo $p$. The $p$-adic order of the nonzero binomial coefficients $\binom{N+M}{N-kp}$ does not depend on $k$:
$$\operatorname{ord}_p \binom{N+M}{N-kp} = -\left\{\frac{N}{p}+\frac{M}{p}\right\} + \left\{\frac{N}{p}-k\right\} + \left\{\frac{M}{p}+k\right\} = -\left\{\frac{N}{p}+\frac{M}{p}\right\} + \left\{\frac{N}{p}\right\} + \left\{\frac{M}{p}\right\}.$$
Recalling that $N=n+m$, $M=2n-m$ with $m=n_1+n_2$ the latter quantity reads
$$\operatorname{ord}_p \binom{3n}{n+n_1+n_2} = \lfloor 3\omega \rfloor - \lfloor \omega+\omega_1+\omega_2 \rfloor - \lfloor 2\omega-\omega_1-\omega_2 \rfloor = 1,$$
where we employed ($\star$) to get
$$\lfloor \omega+\omega_1+\omega_2 \rfloor = \lfloor 2\omega-\omega_1-\omega_2 \rfloor = 0.$$
This means that all binomial coefficients $\binom{N+M}{N-kp}$ are divisible by $p$, thus completing our proof of the divisibility of $\hat{T}(n, n_1, n_2)$ by $p$, and of the lemma. \hfill $\blacksquare$

\textbf{Lemma 3.} Define $\Phi = \Phi_n = \prod_{p \in \mathcal{P}_n} p$ and
$$ L_n = \frac{\operatorname{lcm}(1, 2, \dots, 4n)}{\Phi_n} \in \mathbb{Z}. $$
Then
\begin{equation}
    17 \times L_n \times \frac{34^{-j} A_j}{j} \in \mathbb{Z} \quad \text{for } j \in \{-4n, -4n+1, \dots, 3n\}, j \ne 0.
    \label{eq:lemma3}
\end{equation} 
and $\Phi_n^{-1} \times A_0 \in \mathbb{Z}$.

Asymptotically,
$$ \lim_{n \to \infty} \frac{\log \Phi_n}{n} = \frac{\Gamma'(2/3)}{\Gamma(2/3)} - \frac{\Gamma'(1/2)}{\Gamma(1/2)} = \frac{\pi}{2\sqrt{3}} - \log \frac{3\sqrt{3}}{4} = 0.64527561\dots $$.

\textbf{Proof.} Obviously,
$$ \operatorname{lcm}(1, 2, \dots, 4n) \times \frac{1}{j} \in \mathbb{Z} \quad \text{for } j \in \{-4n, -4n+1, \dots, 3n\}, j \ne 0 $$
Consider
$$ L_n \times \frac{1}{j/p} = \frac{\operatorname{lcm}(1, 2, \dots, 4n)/j}{\Phi_n/p}$$
If $p|j$,  $p \in \mathcal{P}_n$, we just need to prove that for any $q \in \mathcal{P}_n$  other than $p$ cannot divide $j$. Since $ q \in \mathcal{P}_n $, it follows that  $q > \sqrt{10n}$. If both $p $ and $q $ divide $j$ , then $pq $ divides $j $ and consequently $|j| \geq pq > 10n$ . But this contradicts $|j| \leq 4n $. Therefore, a number $j$ can only have one prime factor from $\mathcal{P}_n$.

Therefore, we can have:
$$L_n \times \frac{1}{j/p} \in \mathbb{Z} \quad \text{if } p|j , p \in \mathcal{P}_n. \quad (\star \star)$$

On the other hand, it follows from  Lemma 2 that
$$ \frac{34^{-j}17 A_j}{p} \in \mathbb{Z} \quad \text{if } p|j, p \in \mathcal{P}_n. \quad (\star \star \star) $$

Combining $(\star \star)$ and $(\star \star \star)$ results in claim \eqref{eq:lemma3}. For $j = 0$, Lemma 2 asserts that $p | A_0$ for any prime $p \in \mathcal{P}_n$. Since $\Phi_n$ is the product of these distinct primes, we conclude that $\Phi_n$ divides $A_0$, yielding $\Phi_n^{-1} \times A_0 \in \mathbb{Z}$. 

By definition, $\log \Phi_n = \sum_{p \in \mathcal{P}_n} \log p$. Let $\theta(x) = \sum_{p \le x} \log p$ be the Chebyshev function. The theoretical intervals dictated by the fractional part condition $1/2 \le \{n/p\} < 2/3$ correspond to $p \in \bigcup_{k=0}^{\infty} (\frac{n}{k+2/3}, \frac{n}{k+1/2}]$. 

Notice that our set $\mathcal{P}_n$ enforces a lower bound $p > \sqrt{10n}$. The contribution of the primes omitted by this truncation is strictly bounded by $\theta(\sqrt{10n})$. By  Prime Number Theorem $\theta(x) = \mathcal{O}(x)$, meaning the omitted sum is $\mathcal{O}(\sqrt{n})$. When divided by $n$, this error term vanishes as $n \to \infty$:
$$ \lim_{n \to \infty} \frac{\mathcal{O}(\sqrt{n})}{n} = 0. $$
Therefore, the asymptotic limit is unaffected by the lower bound, and we can safely sum over the full intervals:
$$ \lim_{n \to \infty} \frac{\log \Phi_n}{n} = \lim_{n \to \infty} \frac{1}{n} \sum_{k=0}^{\infty} \left( \theta\left(\frac{n}{k+1/2}\right) - \theta\left(\frac{n}{k+2/3}\right) \right). $$
By the Prime Number Theorem, $\theta(x) \sim x$ as $x \to \infty$, leading to:
$$ \lim_{n \to \infty} \frac{\log \Phi_n}{n} = \sum_{k=0}^{\infty} \left( \frac{1}{k+1/2} - \frac{1}{k+2/3} \right). $$
Recall the series representation of the Digamma function $\psi(z) = \frac{\Gamma'(z)}{\Gamma(z)}$:
$$ \psi(y) - \psi(x) = \sum_{k=0}^{\infty} \left( \frac{1}{k+x} - \frac{1}{k+y} \right). $$
Setting $x = 1/2$ and $y = 2/3$ yields the desired sum. Using Gauss's Digamma theorem, we evaluate $\psi(2/3) - \psi(1/2) = \frac{\pi}{2\sqrt{3}} - \log \frac{3\sqrt{3}}{4} \approx 0.64527561\dots$, which completes the proof. \hfill $\blacksquare$

\textbf{Lemma 4.} \textit{Write the polynomial $P(x) \in \mathbb{Z}[x]$ in the decomposition from the partial-fraction expansion as}
\begin{equation}
P(x) = \sum_{k=0}^{4n-2} B_k (x+1+4i)^k \quad \text{with } B_k \in \mathbb{Z}[i] \text{ for } k=0, 1, \dots, 4n-2.
\end{equation}
\textit{Then the coefficients satisfy}
\begin{equation}
2^{-(3n-3) + 2k}  \times B_k \in \mathbb{Z}[i] \quad \text{for } k=0, 1, \dots, 4n-2.
\end{equation}

\begin{proof}
If $k \ge \lceil 3n/2 \rceil$, the inclusion follows from $B_k \in \mathbb{Z}[i]$ and the fact that the scaling exponent remains non-negative. Therefore, we only need to verify the property for $k < \lceil 3n/2 \rceil$. 

Since the integrand $R(x)$ has a zero of order $2n$ at $x = -1-4i$, we deduce from the partial fraction decomposition that the Taylor coefficients of $P(x)$ must satisfy:
\begin{equation}
B_k = -\frac{1}{k!} \frac{d^k}{dx^k} \sum_{j=0}^{3n} \left( \frac{A_j}{(17+x)^{j+1}} + \frac{A_j}{(17-x)^{j+1}} \right) \Bigg|_{x=-1-4i}
\end{equation}
\begin{equation}
= -\sum_{j=0}^{3n} \binom{j+k}{k} \left( \frac{(-1)^k A_j}{(17+x)^{j+k+1}} + (-1)^{j+1+k} \frac{A_j}{(x-17)^{j+k+1}} \right) \Bigg|_{x=-1-4i}.
\end{equation}

Substituting the values $17+x = 16-4i = 4(4-i)$ and $x-17 = -18-4i = -2(4-i)(2+i)$ at the point $x = -1-4i$:
\begin{equation}
B_k = -\sum_{j=0}^{3n} \binom{j+k}{k} \left( \frac{(-1)^k A_j}{(4(4-i))^{j+k+1}} + \frac{A_j}{(2(4-i)(2+i))^{j+k+1}} \right).
\end{equation}

Recalling from Lemma 1 that $2^{-(3n+2j-1)} 5^{n-j} 17^{1-j} A_j \in \mathbb{Z}[i]$, we can obtain that:

\begin{equation}
2^{-(3n-3) + 2k} \times (4-i) ^{3n+k+1}\times (2+i)^{3n+k+1} \times B_k \in \mathbb{Z}[i].
\end{equation}

As $B_k \in  \mathbb{Z}[i]$ is already established, we have:
\begin{equation}
2^{-(3n-3) + 2k} \times B_k \in \mathbb{Z}[i].
\end{equation}
\end{proof}

\textbf{Lemma 5.} For the polynomial $P(x)$ in the decomposition (3), we have
$$ 2^{-3n}\times L_n \times i \int_{-1-4i}^{-1+4i} P(x) dx \in \mathbb{Z}. $$

\textbf{Proof.} We first compute the integral using representation (21):
$$ i \int_{-1-4i}^{-1+4i} P(x) dx = i \sum_{k=0}^{4n-2} B_k \int_{-1-4i}^{-1+4i} (x+1+4i)^k dx = i \sum_{k=0}^{4n-2} \frac{B_k}{k+1} (8i)^{k+1} = - \sum_{k=0}^{4n-2} \frac{2^{3k+3} B_k}{k+1} i^k. $$

Using the bounds from Lemma 4, the 2-adic valuation of the $k$-th term in the summation is at least $(3k+3) + (3n-2k-3) - v_2(k+1) = 3n + k  - v_2(k+1)$. Since $m - v_2(m) \ge 1$ for any integer $m = k+1 \ge 1$, the overall power of 2 is strictly bounded below by $3n$. Multiplying by the least common multiple to clear the odd components of $k+1$, we obtain:
\begin{equation}
     2^{-3n} \times \operatorname{lcm}(1, 2, \dots, 4n) \times i \int_{-1-4i}^{-1+4i} P(x) dx \in \mathbb{Z}[i]. 
     \label{eq:2}
\end{equation}

On the other hand, we apply representation (20):
$$ i \int_{-1-4i}^{-1+4i} P(x) dx = i \sum_{k=1}^{4n-1} \left( A_{-k} - \sum_{j=0}^{3n} \binom{j+k-1}{j} \frac{A_j}{34^{j+k}} \right) \int_{-1-4i}^{-1+4i} (x+17)^{k-1} dx $$
$$ = \sum_{k=1}^{4n-1} \left( \frac{A_{-k}}{k} - \frac{A_0}{k 34^k} - \sum_{j=1}^{3n} \binom{j+k-1}{j-1} \frac{A_j}{j 34^{j+k}} \right) 2^{2k+1} \sum_{\substack{l=0 \\ l \text{ odd}}}^k \binom{k}{l} 4^{k-l} (-1)^{(l+1)/2} $$
is a rational number satisfying
\begin{equation}
    17\times\frac{\operatorname{lcm}(1, 2, \dots, 4n)}{\Phi_n} \times i \int_{-1-4i}^{-1+4i} P(x) dx \in 34^{-4n} \mathbb{Z} 
    \label{eq:3}
\end{equation}
on the basis of Lemma 3. Finally, since $\Phi_n$ is comprised of primes strictly greater than $17$ (and therefore coprime to $34$), the two inclusions \eqref{eq:2}and \eqref{eq:3} mutually constrain the denominator, seamlessly combining into Lemma 5. \hfill $\blacksquare$

\vspace{0.5cm}

\textbf{Lemma 6.} For the partial-fraction part in \eqref{eq:4} (without the $j=0$ term), we have
$$ 2^{-3n+1} 5^{n} 17^2 L_n \times i \int_{-1-4i}^{-1+4i} \sum_{j=1}^{3n} A_j \left( \frac{1}{(17+x)^{j+1}} + \frac{1}{(17-x)^{j+1}} \right) dx \in \mathbb{Z}. $$

\textbf{Proof.} This rigorously follows from integrating the partial fractions:
$$ i \sum_{j=1}^{3n} A_j \int_{-1-4i}^{-1+4i} \left( \frac{1}{(17+x)^{j+1}} + \frac{1}{(17-x)^{j+1}} \right) dx $$
$$ = i \sum_{j=1}^{3n} \frac{A_j}{j} \left( \frac{1}{(16-4i)^j} - \frac{1}{(16+4i)^j} - \frac{1}{(18+4i)^j} + \frac{1}{(18-4i)^j} \right) $$
$$ = i \sum_{j=1}^{3n} \frac{A_j}{j} \left( \frac{(4+i)^j}{2^{2j} 17^j} - \frac{(4-i)^j}{2^{2j} 17^j} - \frac{(4+i)^j(2-i)^j}{2^j 5^j 17^j} + \frac{(4-i)^j(2+i)^j}{2^j 5^j 17^j} \right) \in \mathbb{Q}. $$
By the bounds explicitly proven in Lemma 1, we have $A_j \in 2^{3n+2j-1} 5^{j-n} 17^{j-1} \mathbb{Z}$. When paired with the sum, the maximal denominators structurally generated are $2^{2j}$, $5^j$, and $17^j$. Reducing the valuation of $A_j$ by these denominators leaves exact residual lower bounds of $2^{3n-1}$, $5^{-n}$, and $17^{-1}$. Clearing these precise fractional residues requires multiplying the expression by $2^{-3n+1}$, $5^{n}$, and $17^2$, while $L_n$ absorbs the fraction's $1/j$ coefficients as guaranteed by Lemma 3. \hfill $\blacksquare$

Lemma 1 and the integrality of $L_n$ imply that $2^{-3n+1} 5^{n} 17^2 L_n \times A_0 \in \mathbb{Z}$; together with the exact calculation of the residual simple poles integration:
$$ \int_{-1-4i}^{-1+4i} \left( \frac{1}{17+x} + \frac{1}{17-x} \right) dx = \left[ \log(17+x) - \log(17-x) \right]_{-1-4i}^{-1+4i} $$
$$ = (\log(16+4i) - \log(18-4i)) - (\log(16-4i) - \log(18+4i)) $$
$$ = \log\left(\frac{16+4i}{16-4i}\right) + \log\left(\frac{18+4i}{18-4i}\right) $$
$$ = \log\left(\frac{4+i}{4-i}\right) + \log\left(\frac{9+2i}{9-2i}\right) $$
By evaluating the imaginary arguments, this simplifies precisely to:
$$ = 2i \arctan\left(\frac{1}{4}\right) + 2i \arctan\left(\frac{2}{9}\right) = 2i \arctan\left(\frac{1}{2}\right) $$
and Lemmas 5, 6 we are thus led to the following statement.

\vspace{0.5cm}

\textbf{Proposition 1.} For the integrals $I_n$ in (1), we have
$$ 2^{-3n+1} 5^{n} 17^2 L_n \times (-i I_n) \in \mathbb{Z} + \mathbb{Z} \arctan(1/2). $$
\textbf{Proposition 2.} The asymptotics of the integrals $I'_n = -i I_n$ and the coefficients $b_n$ in the representation $I'_n = a_n + b_n \arctan(1/2)$ is as follows:
$$ \lim_{n \to \infty} |I'_n|^{1/n} = |N_2| = 0.005442673\dots $$
and
$$ \lim_{n \to \infty} |b_n|^{1/n} = |N_3| = 2621725.228\dots $$
where
$$ N_1 = 0.002126311\dots, \quad N_2 = 0.005442673\dots, \quad N_3 = 2621725.228\dots $$
are the local extrema of the rational function associated with the saddle-point method.

\textbf{Proof.} 
By saddle-point method, we need to 
 compute the saddle points of this function
$$ g(y) = \frac{y(y^2+30y+289)^2}{(y-289)^3} $$
and with the zeros of the logarithmic derivative
$$ \frac{g'(y)}{g(y)} = \frac{1}{y} + \frac{2(2y+30)}{y^2+30y+289} - \frac{3}{y-289} =\frac{2y^3 - 1445y^2 - 26588y - 83521}{y(y^2+30y+289)(y-289)}. 
$$
The zeros of this numerator are precisely:
$$ y_1 = -14.000199\dots, \quad y_2 = -4.028004\dots, \quad y_3 = 740.528203\dots $$
Then $N_j = |g(y_j)|$ for $j=1, 2, 3$. By the standard steepest descent analysis, the original contour of integration can be deformed to pass through the saddle point $y_2$ on the real axis, which strictly dominates the asymptotic behavior of $I_n$. Similarly, the conjugate contour determining the coefficients $b_n$ is dominated by the complementary saddle point $y_3 > 289$.\hfill $\blacksquare$

\textbf{Irrationality measure of $\arctan(1/2)$.} It follows from Propositions 1 and 2 that the forms
$$ I'_n = 2^{-3n+1} 5^{n} 17^2 L_n (-i I_n) = a'_n + b'_n \arctan\left(\frac{1}{2}\right), \quad \text{where } n = 0, 1, 2, \dots, $$
all have integral coefficients $a'_n, b'_n$ and the asymptotics
$$ \limsup_{n \to \infty} \frac{\log |I'_n|}{n} = \log |N_2| - 3 \log 2 + \log 5 + 4 - \frac{\pi}{2\sqrt{3}} + \log \frac{3\sqrt{3}}{4} = -2.328764\dots $$
and
$$ \lim_{n \to \infty} \frac{\log b'_n}{n} = \log |N_3| - 3 \log 2 + \log 5 + 4 - \frac{\pi}{2\sqrt{3}} + \log \frac{3\sqrt{3}}{4} = 17.664067\dots $$ This implies  that the irrationality measure of $\arctan(1/2)$ is bounded above by
$$ 1 + \frac{17.664067\dots}{2.328764\dots} = 8.585166\dots $$

\end{document}